\documentclass[a4paper,twoside]{article}
\usepackage{amsthm,amscd}

\usepackage{amsmath,amssymb,amsfonts}

\usepackage{graphicx}               
\usepackage{color}                  
\RequirePackage[colorlinks=true,citecolor=blue,linkcolor=blue]{hyperref}
\usepackage{hyperref}
\usepackage[absolute]{textpos} 
\usepackage{multicol}
\usepackage{array}

\theoremstyle{remark}

\newtheorem{rem}{Remark}

\usepackage[T1]{fontenc}

\usepackage{xcolor}
\usepackage[french]{babel}

\newcommand{\beq} {\begin{eqnarray*}}
\newcommand{\eeq} {\end{eqnarray*}}

\theoremstyle{plain}
\newtheorem{prop}{Proposition}[section]

\begin{document}





%
%
%
%
%

\title{Remark on the finite-dimensional character  of certain  results of functional statistics\\
Remarque sur le caractère fini dimensionnel de certains résultats de statistique fonctionnelle 
}

\author{ Jean-Marc Aza\"{\i}s \thanks { Universit\'e de Toulouse, IMT, ESP,
  F31062
Toulouse Cedex 9, France. Email:  jean-marc.azais@math.univ-toulouse.fr } \and
Jean-Claude Fort \thanks{Universit\'e Paris Descartes, SPC, MAP5, 45 rue des saints p\`eres, F-75006 Paris;
Email: jean-claude.fort@parisdescartes.fr}}
\maketitle

{\bf Résumé}\\
Cette note montre qu'une hypothèse concernant les probabilités de petites boules, fréquemment utilisée en statistique fonctionnelle, implique que la dimension de l'espace fonctionnel considéré est finie. Un exemple  de processus $L^2$, ne vérifiant pas cette hypothèse, vient compléter ce résultat.

{\bf Abstract}\\
This note shows that some assumption on small balls probability, frequently used in the domain of functional statistics, implies that the considered functional space is of finite dimension.  To complete this result an example of $L^2$ process is given that does not fulfill this assumption.

%

\section{The result}
In several  functional  statistics papers 
 (cf \cite{Bu},\cite{Da2},\cite{Fe-1},\cite{Fe0},\cite{Fe1},\cite{Fe2}) the following hypothesis is used: 
 
{\it(H) Let  $x$ be a point of the space $\cal X$ where a functional variable  $X$ lives. The space $\cal X$    is equipped with a semi distance  and  $B(x,h)$  is the ball with center $x $and radius $h>0$. 
We set $\varphi_x(h)= \mathbb P(X\in B(x,h))$ and we assume :}
$$
\inf_{h\in[0,1]}\int_0^1\varphi_x(ht)dt/\varphi_x(h)\ge\theta_x>0.
$$
where  the parameter $\theta_x$ is locally bounded away from zero.


 The aim of this note  is to prove that $(H)$ implies that $\cal X$  is of finite dimension. \bigskip
 
{\em Proof:} Without loss of generality we can assume that $\theta_x <1/2$.
Let  $F_x(h)=\int_0^h\varphi_x(t)dt$ we have:

$$\frac{1}{h}\int_0^h\varphi_x(t)dt/\varphi_x(h)=\frac{F_x(h)}{hF_x'(h)}\ge \theta_x.$$
By integration we obtain
\begin{equation}\label{un}
F_x(h)/F_x(1)\ge h^{1/\theta_x}.
\end{equation}
Since $\varphi$  is non decreasing  we have
$$
F_x(h)  \leq h \varphi_x(h)
$$
and thus
$$
\varphi_x(h)  \geq h ^{1/\theta_x -1} F_x(1).
$$

 Let $x\in  \cal X$  such that  the parameter $\theta_y$ has positive lower bound   for $y$ in the ball $B(x,h_0)$. This implies  that $\varphi_x(h_0) >0$. By a scaling we can assume for simplicity that  $h_0=1$. For all $y$ in 
 $B(x,1/4)$, and for all $h\in [1/2,1]$, 
 $$
 \varphi_y(h) \geq \varphi_x(1/4) \geq  (\frac 1 4) ^{1/\theta -1} F_x(1),
 $$
 where $\theta$ is the uniform lower bound  for $\theta_y$ for $y$ in $B(x,1)$.
 By  integration
 $$
 F_y(1) \geq 1/2  (\frac 1 4) ^{1/\theta -1} F_x(1)
 $$
 and 
 $$
  \varphi_y(h)  \geq  h^{1/\theta -1}  F_y(1) \geq 1/2  (\frac 1 4) ^{1/\theta -1}  h^{1/\theta -1} F_x(1)$$
   This implies that there exist at most $  \mathcal{O}( h ^{1- \theta } )$ disjoints ball of radius $h$ in $B(x,1/4)$. 
    Then the same set of balls but with radius 2r is a covering, which  implies in turn that the box (or entropy) dimension of B(x,1/4) is finite. Since the Hausdorff dimension is smaller than the box dimension it is also finite. \bigskip

\begin{rem}
 Suppose, in addition,  that  the probability distribution function  of $X$  satisfying $H$  admits a density  with respect  to some natural  positive measure $\mu$  and suppose that this density is locally  upper- and lower -bounded by $M$ and $m$ respectively :
 $$ m\mu(B(x,h))\le\varphi_x(h)\le M\mu(B(x,h)).
 $$ 
 We have then  \begin{equation*}
 \ \exists d>0, C_2>0,h_0>0, \ \ 0\le h\le h_0 \Longrightarrow  1/C_2 h^d \le \mu(B(x,h))\le C_2h^d.
 \end{equation*}
\end{rem}
which means that $\mu$ shares the same property as the distribution of $X$. As for example this property does not hold for the Wiener measure. \bigskip 

 We may note that  the classical  random processes:  Brownian motion and  more general Gaussian processes for which the probability of small balls is known do not satisfy {\it H} , see \cite{Lishao}  for more details.
 
  \section{Random series of functions }
 
 Here we consider $ {X}(t)$ a stochastic process constructed  as a random series of functions. More precisely
 \begin{equation}\label{f:2}
  {X}(t)=\sum_{n=1}^\infty \alpha_n Z_n\varphi_n(t)
  \end{equation}
   where
   \begin{itemize}
  \item $(Z_n,n\ge 1)$ is a sequence of independent real random variables, with mean $0$ and variance $1$, with absolutely continuous   densities $(f_n, n\ge 1)$ $w.r.t$ to the Lebesgue measure.  This is the case for most of the usual continuous distributions of probability: exponential, normal, polynomial, gamma, beta etc...
  \item  $(1;\varphi_n, n\ge 1)$ is an orthonormal basis of $L^2([0,1])$.
 \item $(\alpha_n, n\ge 1)$ is a sequence of positive real numbers  $\sum_{n=1}^{+\infty} \alpha_n^2 <\infty$.
 \end{itemize}
 
 The sum (\ref{f:2}) converges in $L^2([0,1])$. 
 
In particular the form  (\ref{f:2}) covers all  the Gaussian processes  through the Karhunen-Lo\`eve  decomposition. 
 
 Then we have : 
 
 \begin{prop}
 $
 \displaystyle \lim_{h\to 0} h^d \mathbb P(\|{ X}\|_2 \le h)=0$ for any $d\ge 0$. So that the process ${X}(t)$ does not fulfill the assumption {\it(H)} for the $L^2$ norm  at the point zero.
 \end{prop}
 The proof is based on the properties of the convolution:
 
 Let $f$ and $g$ the  absolutely continuous densities  of probability of two independent random variables $U$ and $V$. Then $U^2$ (resp. $V^2$) has  the density
 $$
   p_{U^2} (u) = \frac{\widetilde{f}(\sqrt{u})}{ \sqrt u} \ \ \ \mbox{ resp. }  p_{V^2} (v) = \frac{\widetilde{g}(\sqrt{v})}{ \sqrt v},
 $$
 where $\widetilde f$ and $\widetilde g$ are the symmetrized of $f$ and $g$,
 $\widetilde f = 1/2(f(x) +f(-x))$. It follows that $U^2+V^2$ has for density
 $$
  C(x)=\int_0^x             \frac{\widetilde{f}(\sqrt{u})}{ \sqrt u}      \frac{\widetilde{g}(\sqrt{x-u})}{ \sqrt{x-u}}   du = \int_0^1 \frac{1}{\sqrt{v(1-v)}}\tilde g(\sqrt{(1-v)x})\tilde f(\sqrt{vx})dv ,
  $$
  for  $ x\ge 0$ and $ 0$ elsewhere.

  It is easy to see that for any $0<A<B$, the function $C$ is Lipschitz on [A,B]. At $x=0$ it takes the  value $0$ with  right limit $\beta(\frac{1}{2},\frac{1}{2}) f(0)g(0)$. 
Now if we make the convolution product of two such functions $C_1$ and $C_2$, we obtain a function vanishing for $x\le 0$, continuous at $0$ and Lipschitz on any compact interval of $\mathbb R$, thus absolutely continuous. Using a classical result, making the convolution product of $k$ such absolutely continuous functions yields a ${\cal C}^{k-1}$ function whose $(k-1)^{\mbox{th}}$ derivative is absolutely continuous.

Then, applying iteratively this result we conclude that the density of the variable $\displaystyle \|{ X}\|_2^2=\sum_{n=1}^\infty \alpha_n^2 Z_n^2$ is infinitely derivable with all its derivatives null at $0$. The claimed property follows. 

\begin{rem} If the process $ X$ lives in $L^p[0,1], p\in(2,\infty]$, we have
 $\|{X}\|_2^2 \le\|{X}\|_p ^2$, thus:
 $$ \mathbb P(\|{ X}\|_p \le h)\le\mathbb P(\|{ X}\|_2\le h).$$  
\end{rem}


%

\renewcommand{\refname}{References}

\end{document}